\documentclass[12pt]{amsart}
\usepackage{amsmath,amsfonts,amstext,amsthm,epsfig,xcolor}

\newtheorem{thm}{Theorem}[section]

\newtheorem{prop}[thm]{Proposition}
\newtheorem{hyp}[thm]{Hypothesis}
\newtheorem{rem}[thm]{Remark}
\numberwithin{equation}{section}

\newcommand{\umlaut}{\"}

\DeclareMathOperator{\Lip}{Lip}
\DeclareMathOperator{\Log}{Log}

\begin{document}

\title[Reconstruction of the initial condition]
{Reconstruction of the initial condition in parabolic equations with Log-Lipschitz coefficients}%

\author{Daniele Del Santo and Martino Prizzi}

\address{D. Del Santo, Dipartimento di Matematica, Informatica e Geoscienze, Universit\`a degli Studi di Trieste, Via A. Valerio 12 - 34100 Trieste, Italy} \email{delsanto@units.it}
\address{M. Prizzi, Dipartimento di Matematica, Informatica e Geoscienze, Universit\`a degli Studi di Trieste, Via A. Valerio 12 - 34100 Trieste, Italy} \email{mprizzi@units.it}

\begin{abstract}
We consider a parabolic equation whose coefficients are Log-Lipschitz continuous in $t$ and Lipschitz continuous in $x$.
Combining a recent conditional stability result with a well posed variational problem, we reconstruct the initial condition of an unknown solution from a rough measurement at the final time.
\end{abstract}
\maketitle
\section{Introduction}
In this paper we consider the strictly parabolic equation 
\begin{equation} \label{intro1}
 \partial_t u -\sum\limits_{j,k=1}^n \partial_{x_j}(a_{jk}(t,x)\partial_{x_k}u) = 0
\end{equation} 
on the strip $[0,T] \times \Bbb R^n$. We assume that the coefficients $a_{ij}$'s are Lipschitz continuous with respect to $x$ and Log-Lipschitz continuous with respect to $t$ (we shall make these assumptions more precise in the next section). It is well known that the backward Cauchy problem for (\ref{intro1}) is {\em ill posed} (see e.g. \cite{Isakov} and the references therein).
Yet, our aim is to reconstruct the initial condition $\bar u(0)$ of an unknown solution $\bar u(t)$ of (\ref{intro1}), provided we have an {a priori} bound on $\bar u(0)$ and we can measure with arbitrary accuracy its final configuration $\bar u(T)$.  More precisely, suppose that for every $\delta>0$ we can perform a measurement $g$ of $\bar u(T)$ such that 
\begin{equation}\label{intro2}
\|g-\bar u(T)\|_{L^2}\leq\delta.
\end{equation}
Moreover, suppose that we know {\em a priori} that $\|\bar u(0)\|_{H^1}\leq E$ for some $E>0$. We are interested in finding a {\em computable} approximation of $\bar u(0)$. By {\em computable} we mean that it is obtained by solving some {\em well posed} problem. There is a vast literature on the subject when the coefficients $a_{ij}(t,x)$ are Lipschitz continuous in $t$. The procedure consists typically in two steps: 
\begin{itemize}
\item first, one proves some kind of {\em well behaveness} or {\em conditional stability} for backward solutions of (\ref{intro1}). 
The concept of {\it conditional stability} is central in inverse problem theory and was introduced by Tikhonov in 1943 in the paper \cite{Tikhonov}, while the more stringent concept of {\em well behaveness} was introduced by John in 1960 in the paper \cite{John1960}. Roughly, conditional stability means continuous dependence on data, while well behaveness means H\umlaut older continuous dependence on data, subject in both cases to some a priori bounds. We refer to \cite{CDSP2022} for a (partial) history and a discussion of these concepts.
\item second, one searches for a technique to construct some solution of (\ref{intro1}) to which the conditional stability result applies.
\end{itemize}

The first step is usually achieved by exploiting the {\em logarithmic convexity technique}, developed by Agmon and Nirenberg in \cite{AN1963}. The Lipschitz continuity of the coefficients $a_{ij}(t,x)$ with respect to $t$ plays a crucial role in the computations involved in the logatithmic convexity technique. For the LogLipschitz case, in the paper
\cite{CDSP2022} the following {\em conditional stability} 
result was proved through weighted energy-type estimates obtained by mean of microlocal analysis tools:

{\it for every $E>0$ there exist $\rho>0$,
$0<\theta\leq1$ and $K>0$ such that, if $u$ is a
solution of (\ref{intro1}) on $[0,T]$ with $\|u(0,\cdot)\|_{H^1}\leq E$ and
$\|u(T,\cdot)\|_{L^2}\leq\rho$, then 
\begin{equation}\label{intro_dis}\sup_{t\in[0,T]}\|u(t,\cdot)\|_{L^2}\leq
K\frac1{|\log\|u(T,\cdot)\|_{L^2}|^\theta}.\end{equation}}

Concerning the second step, if it were possible to solve equation (\ref{intro1}) backward in time with final condition $u(T)=g$, then by (\ref{intro_dis}) we would get that $u(0)$ is closed to $\bar u(0)$, provided $\|u(0)\|_{H^1}\leq E$ and $g$ is sufficiently closed to $\bar u(T)$. However, equation (\ref{intro1}) with final condition $u(T)=g$ in general has no solution, due to the regularizing effect of equation (\ref{intro1}) forward in time, and to the fact that $g$ does not possess any regularity, since it is the output of a measurement. There are various strategies to overcome this obstruction. We mention the technique of quasi reversibility (see e.g. \cite{Ew}), which consists in perturbing the equation to make it solvable backward in time, and the technique of Fourier truncation, which consists in approximating $g$ with a very regular function obtained truncating its Fourier transform. The latter technique was exploited in \cite{CDSP2022} to illustrate the conditional stability result  through an example inspired by \cite{FXQ}. In that example the coefficients $a_{ij}$'s depended only on the time variable $t$. 

On the other hand, when the coefficients $a_{ij}$'s depend on $x$ but are independent of $t$, another technique which is frequently used is the so called {\em nonlocal quasi-boundary value method}, consisting in solving the {\em well posed} nonlocal problem
\begin{equation}
\begin{cases}
 \partial_t u -\sum\limits_{j,k=1}^n \partial_{x_j}(a_{jk}(x)\partial_{x_k}u) = 0\\
 \alpha u(0)+u(T)=g
\end{cases}
\end{equation}
where $\alpha$ is a small parameter to be chosen in a suitable way, depending on the parameter $\delta$ of (\ref{intro2}). This technique was firstly developed by Showalter in \cite{showalter}. See also \cite{clark, denche, HDL, HDS, liu}. In \cite{haovanduc} H\`ao and Duc proposed a version of the quasi-boundary value method suited for non autonomous equations. The estimates obtained by H\`ao and Duc however cannot be applied to the case of Log-Lipschitz coefficients, since in their proof the logarithmic-convexity estimates for the solutions play a crucial role. Moreover, the result of H\`ao and Duc furnishes an approximation of $\bar u(t)$ for all $t\in ]0,T]$, but gives no information at $t=0$.  

In the present paper we further develope some ideas of H\`ao and Duc and obtain a full approximation of $\bar u(t)$ up to $t=0$ when the coefficients $a_{ij}(t,x)$ are Log-Lipschitz continuous with respect to $t$. To be more precise, if $\bar u(\cdot)$ is an unknown solution of (\ref{intro1}) satisfying a $H^1$ a-priori bound, and $\delta$ and $g$ are given as in (\ref{intro2}), we construct a solution $u(t)$ of (\ref{intro1}) such that
\begin{equation}\label{intro_rate}
\sup_{t\in[0,T]}\|u(t)-\bar u(t)\|_{L^2}\leq K'\frac1{|\log \delta|^\theta},
\end{equation}
where $\theta$ is the same exponent of (\ref{intro_dis}). 
In Section 4 we shall discuss briefly how this convergence rate relates to similar results concerning autonomous equations and non-autonomous equations whose coefficients are Lipschitz continuous in $t$.

\section{Preliminaries} \label{sec:prel}

We consider the parabolic equation 
\begin{equation} \label{pareq}
 \partial_t u -\sum\limits_{j,k=1}^n \partial_{x_j}(a_{jk}(t,x)\partial_{x_k}u) = 0
\end{equation} 
on the strip $[0,T] \times \Bbb R^n$. 
\begin{hyp} \label{mainhypothesis} We assume throughout the paper that:
\begin{enumerate}
\item for all $(t,x) \in [0,T] \times \Bbb R^n$ and for all $j,k = 1, \dots, n$, \begin{equation*}
a_{jk}(t,x) = a_{kj}(t,x);
\end{equation*}
\item there exists $\kappa$, $0<\kappa\leq1$, such that for all $(t,x,\xi) \in [0,T] \times \Bbb R^n \times \Bbb R^n$, \begin{equation}
\kappa |\xi|^2 \leq \sum\limits_{j,k=1}^n a_{jk}(t,x)\xi_j \xi_k \leq \frac{1}{\kappa} |\xi|^2;
\end{equation}
\item for all $j,k=1,\dots,n$, $a_{jk} \in \Log\Lip([0,T],L^\infty(\Bbb R^n)) \cap L^\infty([0,T], \Lip(\Bbb R^n))$. 
\end{enumerate}
We set \begin{eqnarray*}
&& A_{LL} := \sup\Big\{ \frac{|a_{jk}(t,x)-a_{jk}(s,x)|}{|t-s|(1+|\log|t-s||)} \mid j,k = 1, \dots, n, \\
&& \qquad \qquad \qquad\qquad\qquad t,s \in [0,T],\, x\in  \Bbb R^n, \, 0 < |s-t| \leq 1 \Big\}, \\[0.3 cm]
&& A := \sup\{ \|\partial_x^\alpha a_{jk}(t,\cdot)\|_{L^\infty} \mid |\alpha| \leq 1, \, t \in [0,T] \}.
\end{eqnarray*}\end{hyp}

Throughout the paper, whenever $X$ is a Hilbert space we denote by $\langle u|v\rangle_X$ the scalar product of $u$ and $v\in X$ and by $\|\cdot\|_X$ the corresponding norm. Moreover, we denote by $\langle\phi|u\rangle_{(X^*,X)}$ the action of the continuous linear form $\phi\in X^*$ on $u\in X$.

For each $t\in[0,T]$ we define the symmetric, non-negative bilinear form
\begin{equation}
a(t;u,v):=\int_{\Bbb R^n}\sum\limits_{j,k=1}^n a_{jk}(t,x)\partial_{x_j}u(x)\,\partial_{x_k}v(x)\,dx, \quad u,v\in H^1(\Bbb R^n)
\end{equation}
and we notice that for all $u\in H^1(\Bbb R^n)$ we have
\begin{equation}
|a(t;u,v)|\leq\frac1\kappa\|u\|_{H^1}\|v\|_{H^1}
\end{equation}
and
\begin{equation}
\kappa\|\nabla u\|_{L^2}^2\leq a(t;u,u)\leq\frac1\kappa\|\nabla u\|_{L^2}^2,
\end{equation}
whence 
\begin{equation}
\kappa\|u\|_{H^1}^2\leq a(t;u,u)+\|u\|_{L^2}^2\leq\frac1\kappa\|u\|_{H^1}^2,\quad u\in H^1(\Bbb R^n).
\end{equation}
For each $t\in[0,T]$ the bilinear form $a(t;\cdot,\cdot)$ defines a bounded linear operator $\mathcal A(t):H^{1}(\Bbb R^n)\to H^{-1}(\Bbb R^n)$ by the formula
\begin{equation}
\langle \mathcal A(t)u,v\rangle_{(H^{-1},H^1)}:=a(t;u,v), \quad u,v\in H^1(\Bbb R^n).
\end{equation}
We denote by $A(t)$ the part of $\mathcal A(t)$ in $L^2(\Bbb R^n)$.
By regularity theory for elliptic partial differential equations (see e.g. \cite[Thms. 8.8 and 8.12]{GilTrud}), for each $t\in[0,T]$ the operator $A(t)$ is self-adjoint and non negative in $L^2(\Bbb R^n)$, with domain $H^2(\Bbb R^n)$, and its explicit expression is 
\begin{equation}
A(t) u=-\sum\limits_{j,k=1}^n \partial_{x_j}(a_{jk}(t,x)\partial_{x_k}u).
\end{equation}
The dependence on $t$ of the operator $A(t)$ is better than H\umlaut older continuous, so one can apply the abstract theory of linear parabolic equations  and obtain well posedness of the initial value Cauchy problem associated with (\ref{pareq})
in the space $H^\theta(\Bbb R^n)$ for every $0\leq\theta\leq2$. More precisely, we have the following theorem (see e.g. \cite[Thm. 4.4.1]{Amann}; see also \cite{Tanabe} and \cite{Yagi}):

\begin{thm}\label{sol-op} There exists a family of {\em evolution operators} $U(p,s)$, $0\leq s\leq p\leq T$, for equation (\ref{pareq}), satisfying the following properties:
\begin{enumerate}
\item for every $s,p$ with $0\leq s\leq p\leq T$ and for every $\theta$ with $0\leq\theta\leq2$, $U(p,s): H^\theta(\Bbb R^n)\to H^\theta(\Bbb R^n)$ is a bounded linear operator, and there exists $M>0$ (independent of $p,s$ and $\theta$) such that $\|U(p,s)\|_{\mathcal L(H^\theta,H^\theta)}\leq M$;
\item for every $p,s$ with $0\leq s< p\leq T$ and for every $\theta_1,\theta_2$ with $0\leq\theta_1<\theta_2\leq2$, $U(p,s): H^{\theta_1}(\Bbb R^n)\to H^{\theta_2}(\Bbb R^n)$ is a bounded linear operator, and there exists $M>0$ (independent of $p,s$ and $\theta_1,\theta_2$) such that $\|U(p,s)\|_{\mathcal L(H^{\theta_1},H^{\theta_2})}\leq M(1+(p-s)^{-(\theta_2-\theta_1)})$;
\item for every $s\in[0,T]$, $U(s,s)=I$ and
for every $s,q,p$ with $\leq s\leq q\leq p\leq T$, $U(p,q)U(q,s)=U(p,s)$;
\item for every $\theta\in[0,2]$, $s\in[0,T[$ and $u_s\in H^\theta(\Bbb R^n)$, the function $u(p):=U(p,s)u_s$ is continuous from $[s,T]$ to $H^\theta(\Bbb R^n)$, continuous from $]s,T]$ to $H^2(\Bbb R^n)$, differentiable from $]s,T]$ to $L^2(\Bbb R^n)$, and satisfies (\ref{pareq}) in $L^2(\Bbb R^n)$ on $]s,T]$ with $u(s)=u_s$. If $\theta=2$, $u(p)$ is differentiable in $L^2(\Bbb R^n)$ up to $p=s$, and if $\theta\geq 1$ is differentiable at least in $H^{-1}(\Bbb R^n)$ up to $p=s$.
\end{enumerate}\qed
\end{thm}

Together with (\ref{pareq}), we consider the adjoint problem 
\begin{equation} \label{adjpareq}
 \partial_t v +\sum\limits_{j,k=1}^n \partial_{x_j}(a_{jk}(t,x)\partial_{x_k}v) = 0
\end{equation} 
on the strip $[0,T] \times \Bbb R^n$.  This problem is well posed backward in time and it defines an evolution family $V(s,p)$, $0\leq s\leq p\leq T$, where $v(s):=V(s,p)v_p$ is the solution of (\ref{adjpareq}) on $[0,p]$ with $v(p)=v_p$, for every $v_p\in L^2(\Bbb R^n)$. The family $V(s,p)$ satisfies the same properties of $U(p,s)$ with reversed time direction.

\begin{prop}For $0\leq s\leq p\leq T$ the following relation holds:
\begin{equation}V(s,p)=U(p,s)^*,\end{equation}
where $U(p,s)^*$ is the adjoint of $U(p,s)$.\end{prop}
\begin{proof}
Let $u_s,v_p\in L^2(\Bbb R^n)$ and set $u(t):=U(t,s)u_s$, $v(t):=V(t,p)v_p$, $t\in[s,p]$. For $t\in]s,p[$ we have:
\begin{multline*}
\frac{d}{dt}\langle u(t),v(t)\rangle_{L^2}=\langle \partial_t u(t),v(t)\rangle_{L^2}+\langle u(t),\partial_t v(t)\rangle_{L^2}\\
=\langle -A(t)u(t),v(t)\rangle_{L^2}+\langle u(t),A(t)v(t)\rangle_{L^2}\\=a(t;u(t),v(t))-a(t;u(t),v(t))=0.
\end{multline*}
Therefore the map $t\mapsto \langle u(t),v(t)\rangle_{L^2}$ is constant in $[s,p]$. It follows that $\langle u(p),v(p)\rangle_{L^2}=\langle u(s),v(s)\rangle_{L^2}$, that is
$$
\langle U(p,s) u_s,v_p\rangle_{L^2}=\langle u_s,V(s,p)v_p\rangle_{L^2}.
$$
Since $u_s$ and $v_p$ are arbitrary, this means that $V(s,p)=U(p,s)^*$.
\end{proof}

\begin{prop}\label{lax-mil}
Let $\alpha>0$; for $u,v\in H^1(\Bbb R^n)$, let
\begin{equation}\label{bil-form}
b_\alpha(u,v):=\alpha \left(a(0;u,v)+\langle u,v\rangle_{L^2}\right)+\langle U(T,0)u,U(T,0)v\rangle_{L^2}.
\end{equation}
Then:
\begin{enumerate}
\item for all $\phi\in H^{-1}(\Bbb R^n)$ there exists a unique $u_{\alpha,\phi}\in H^1(\Bbb R^n)$ such that
\begin{equation*}b_\alpha(u_{\alpha,\phi},h)=\langle\phi,h\rangle_{(H^{-1},H^1)}\quad\text{for all $h\in H^1(\Bbb R^n)$}\end{equation*}
\item $\|u_{\alpha,\phi}\|_{H^{1}}\leq\frac1{\kappa\alpha}\|\phi\|_{H^{-1}}$
\item $u_{\alpha,\phi}$ achieves
\begin{equation*}
\min_{u\in H^1}\left\{\frac12 b_\alpha(u,u)-\langle\phi,u\rangle_{(H^{-1},H^1)}\right\}
\end{equation*}
\end{enumerate}
\end{prop}
\begin{proof} We have that 
\begin{equation}\label{e1}|b_\alpha(u,v)|\leq\frac\alpha\kappa\|u\|_{H^1}\|v\|_{H^1}+M\|u\|_{L^2}\|v\|_{L^2}\end{equation}
and
\begin{equation}\label{e2}b_\alpha(u,u)\geq\alpha\kappa\|u\|_{H^1}^2.\end{equation}
The thesis follows by the Lax-Milgram theorem (see e.g. \cite[p. 138]{brezis}).\end{proof}

\section{The main result} \label{sec:result}

Let $\bar u\colon[0,T]\to H^1(\Bbb R^n)$ be an unknown solution of (\ref{pareq}). Suppose that we know {\it a priori} that 
$\|\bar u(0)\|_{H^1}\leq E$. Suppose also that for all $\delta>0$ we can measure $\bar u(T)$ with an error smaller than $\delta$, that is we can perform a measurement $g$ of $\bar u(T)$ in such a way that $\|g-\bar u(T)\|_{L^2}\leq\delta$. Our goal is to reconstruct $\bar u(0)$ from $g$ with an error depending on $\delta$ and tending  to $0$ as $\delta\to 0$. Our main tool is the following conditional stability result {\em up to $0$}, that we proved in \cite{CDSP2022}.

\begin{thm}[{\cite[cf Thm. 5.3]{CDSP2022}}]\label{res_jde}
Assume Hypothesis \ref{mainhypothesis} is satisfied. Let $E>0$. There exist positive constants $\rho_E$, $\theta_E$ and $K_E$, depending only on $A_{LL}$, $A$, $\kappa$, $T$ and $E$, such that if $u,v \in C^0([0,T],H^1)\cap C^1(]0,T],L^2)$ are solutions of \eqref{pareq} satisfying $\|u(0,\cdot)\|_{H^1}\leq E$, $\|v(0,\cdot)\|_{H^1}\leq E$ and $\|u(T,\cdot)-v(T,\cdot))\|_{L^2} \leq \rho_E$, then the inequality \begin{eqnarray*}
\sup_{t\in[0,T]}\|u(t,\cdot)- v(t,\cdot)\|_{L^2} \leq K_E \frac1{|\log\|u(T,\cdot)-v(T,\cdot)\|_{L^2}|^{\theta_E}}
\end{eqnarray*} holds true. \qed\end{thm}

If we settle with reconstructing $u(t)$ {\em for $t$ away from $0$}, a better estimate (intermediate between H\umlaut older and logarithmic) is available:
\begin{thm}[{\cite[cf Thm. 2.7]{CDSP2022}}]\label{res_jde_bis}
Assume Hypothesis \ref{mainhypothesis} is satisfied. Let $E>0$ and $t_*\in]0,T[$. There exist positive constants $\rho_{E,t*}$, $\mu_{E,t*}$, $K_{E,t*}$ and $N_{E,t*}$, depending only on $A_{LL}$, $A$, $\kappa$, $T$, $t_*$ and $E$, such that if $u,v \in C^0([0,T],H^1)\cap C^1(]0,T],L^2)$ are solutions of \eqref{pareq} satisfying $\|u(0,\cdot)\|_{L^2}\leq E$, $\|v(0,\cdot)\|_{L^2}\leq E$ and $\|u(T,\cdot)-v(T,\cdot))\|_{L^2} \leq \rho_{E,t_*}$, then the inequality \begin{multline*}
\sup_{t\in[t_*,T]}\|u(t,\cdot)- v(t,\cdot)\|_{L^2} \\\leq K_{E,t_*} \exp\left({-N_{E,t_*}|\log\|u(T,\cdot)-v(T,\cdot)\|_{L^2}|^{\mu_{E,t_*}}}\right)
\end{multline*} holds true. \qed\end{thm}

\begin{rem} Unlike in the Lipschitz case, in the Log-Lipschitz case it is impossible to obtain a conditional stability estimate of H\umlaut older type for $t\geq t_*>0$ (a counterexample was explicitely constructed in \cite{DSP2009}).\end{rem}

In order to apply Theorems \ref{res_jde} and \ref{res_jde_bis} we need to construct a solution $u(\cdot)$ of (\ref{pareq}) such that:
\begin{itemize}
\item $\|u(0)-\bar u(0)\|_{H^1}$ is controlled in terms of $E$;
\item $\|u(T)-\bar u(T)\|_{L^2}$ is small and controlled in terms of $\delta$.
\end{itemize}

\begin{thm}\label{appr} Let $\bar u\colon[0,T]\to H^1(\Bbb R^n)$ be a solution of (\ref{pareq}) with 
$\|\bar u(0)\|_{H^1}\leq E$. Let $\delta>0$ and let $g\in L^2(\Bbb R^n)$ be such that $\|g-\bar u(T)\|_{L^2}\leq\delta$. Then there exists $u_{g,\delta}\in H^1(\Bbb R^n)$ such that, setting $u(t):= U(t,0)u_{g,\delta}$, the following estimates are satisfied:
\begin{enumerate}
\item $\|u(0)-\bar u(0)\|_{H^1}\leq\frac1{\kappa^2}(1+E)$;
\item $\|u(T)-\bar u(T)\|_{L^2}\leq\frac{M^{1/2}}{\kappa}(1+E)\delta^{1/2}$.
\end{enumerate}
The element $u_{g,\delta}$ is obtained as the unique solution of the {\em well posed} problem
\begin{equation*}
b_\alpha(u,h)=\langle u(T,0)^*g,h\rangle_{L^2}\quad\text{for all $h\in H^1(\Bbb R^n)$}
\end{equation*}
with $\alpha=\delta M\kappa$, where $b_\alpha(\cdot,\cdot)$ is the bilinear form defined by (\ref{bil-form}) and $M$ is the constant in Theorem \ref{sol-op}.
\end{thm}
\begin{proof}
Define $\phi\in H^{-1}(\Bbb R^n)$ as follows:
\begin{equation}
\langle\phi,h\rangle_{(H^{-1},H^1)}:=\langle U(T,0)^*g,h\rangle_{L^2}
\end{equation}
Then by Proposition \ref{lax-mil} for all $\alpha>0$ there exists a unique $u_{g,\alpha}\in H^1(\Bbb R^n)$ such that
\begin{equation*}
b_\alpha(u_{g,\alpha},h)=\langle U(T,0)^*g,h\rangle_{L^2}\quad\text{for all $h\in H^1(\Bbb R^n)$}, 
\end{equation*}
that is
\begin{multline}\label{ingred1}
\alpha \left(a(0;u_{g,\alpha},h)+\langle u_{g,\alpha},h\rangle_{L^2}\right)+\langle U(T,0)u_{g,\alpha},U(T,0)h\rangle_{L^2}\\=\langle U(T,0)^*g,h\rangle_{L^2}\quad\text{for all $h\in H^1(\Bbb R^n)$}. 
\end{multline}
Define $u(t):=U(t,0)u_{g,\alpha}$ and $z_0:=u_{g,\alpha}-\bar u(0)=u(0)-\bar u(0)$. Set $g_1:=U(T,0)^*\bar u(T)=U(T,0)^*\,U(T,0)\bar u(0)$ and $g_2:=U(T,0)^*g$. Now we trivially observe that 
\begin{multline}\label{ingred2}
\alpha \left(a(0;\bar u(0),h)+\langle \bar u(0),h\rangle_{L^2}\right)+\langle U(T,0)\bar u(0),U(T,0)h\rangle_{L^2}\\=\alpha \left(a(0;\bar u(0),h)+\langle \bar u(0),h\rangle_{L^2}\right)+\langle g_1,h\rangle_{L^2}\quad\text{for all $h\in H^1(\Bbb R^n)$}. 
\end{multline}
Subtracting (\ref{ingred2}) from (\ref{ingred1}), we get
\begin{multline}\label{ingred3}
\alpha \left(a(0;z_0,h)+\langle z_0,h\rangle_{L^2}\right)+\langle U(T,0)z_0,U(T,0)h\rangle_{L^2}\\=-\alpha \left(a(0;\bar u(0),h)+\langle \bar u(0),h\rangle_{L^2}\right)+\langle g_2-g_1,h\rangle_{L^2}\quad\text{for all $h\in H^1(\Bbb R^n)$}. 
\end{multline}
Choosing $h=z_0$ in (\ref{ingred3}) and taking into account the estimates (\ref{e1}) and (\ref{e2}), we get
\begin{multline}\label{xy}
\alpha\kappa\|z_0\|_{H^1}^2+\|U(T,0)z_0\|_{L^2}^2\\\leq\frac\alpha\kappa\|\bar u(0)\|_{H^1}\|z_0\|_{H^1}+
\|g_2-g_1\|_{L^2}\|z_0\|_{L^2}.
\end{multline}
It follows that
\begin{multline*}
\alpha\kappa\|z_0\|_{H^1}\leq \|g_2-g_1\|_{L^2}+\frac\alpha\kappa\|\bar u(0)\|_{H^1}\\
=\|U(T,0)^*(\bar u(T)-g)\|_{L^2}+\frac\alpha\kappa\|\bar u(0)\|_{H^1}\\
\leq M\|\bar u(T)-g\|_{L^2}+\frac\alpha\kappa\|\bar u(0)\|_{H^1}\leq M\delta+\frac\alpha\kappa E.
\end{multline*}
Therefore we have
\begin{equation}
\|z_0\|_{H^1}\leq\frac{M\delta}{\alpha\kappa}+\frac{E}{\kappa^2}.
\end{equation}
Now we choose $\alpha=\alpha(\delta)=M\kappa\delta$. With this choice we finnaly get
\begin{equation}
\|z_0\|_{H^1}=\|u(0)-\bar u(0)\|\leq\frac1{\kappa^2}(1+E)
\end{equation}
On the other hand, by (\ref{xy}) and with $\alpha=M\kappa\delta$, we have
\begin{multline*}
\|u(T)-\bar u(T)\|_{L^2}^2=\|U(T,0)z_0\|_{L^2}^2\\
\leq\frac\alpha\kappa\|\bar u(0)\|_{H^1}\|z_0\|_{H^1}+\|g_2-g_1\|_{L^2}\|z_0\|_{L^2}\\
\leq\frac\alpha\kappa E\frac1{\kappa^2}(1+E)+M\delta\frac1{\kappa^2}(1+E)\\
\leq M\delta\frac{E}{\kappa^2}(1+E)+M\delta\frac1{\kappa^2}(1+E)=\frac{M}{\kappa^2}(1+E)^2\delta,
\end{multline*}
whence
\begin{equation}
\|u(T)-\bar u(T)\|_{L^2}\leq\frac{M^{1/2}}{\kappa}(1+E)\delta^{1/2}
\end{equation}
The thesis follows defining $u_{g,\delta}:=u_{g,\alpha(\delta)}$.
\end{proof}

Combining theorems \ref{appr} and \ref{res_jde} we obtain the desired reconstruction result:

\begin{thm}\label{rec_th}
Let $\bar u\colon[0,T]\to H^1(\Bbb R^n)$ be an unknown solution of (\ref{pareq}) with $\|\bar u(0)\|_{H^1}\leq E$. Let $E':=(1+E)/\kappa^2$. Let $\rho$, $\theta$ and $K$ be the constants given by Theorem \ref{res_jde} when $E$ is replaced by $E'$. Let $\delta>0$, let $g\in L^2(\Bbb R^n)$ be a measurement of $\bar u(T)$ with $\|g-\bar u(T)\|_{L^2}<\delta$ and let $u_{g,\delta}$ be the element given by Theorem \ref{appr}. Let $u(t):=U(t,0)u_{g,\delta}$. Then, if $\delta$ is so small that $(M^{1/2}/\kappa)(1+E)\delta^{1/2}<\rho$, the estimate
\begin{equation}
\sup_{t\in[0,T]}\|u(t)-\bar u(t)\|_{L^2}\leq K'\frac1{|\log \delta|^\theta}
\end{equation}
holds true (where $K'$ is a constant which can be explicitly expressed in terms of $K$, $\theta$, $\kappa$, $M$ and $E$).\qed
\end{thm}

Combining theorems \ref{appr} and \ref{res_jde_bis} we obtain the a better approximation rate for $t$ away from $0$:

\begin{thm}\label{rec_th_bis}
Let $\bar u\colon[0,T]\to H^1(\Bbb R^n)$ be an unknown solution of (\ref{pareq}) with $\|\bar u(0)\|_{H^1}\leq E$. Let $t*\in]0,T[$ and let $E':=(1+E)/\kappa^2$. Let $\rho$, $\mu$, $K$ and $N$ be the constants (depending on $E'$ and $t_*$) given by Theorem \ref{res_jde_bis} when $E$ is replaced by $E'$. Let $\delta>0$, let $g\in L^2(\Bbb R^n)$ be a measurement of $\bar u(T)$ with $\|g-\bar u(T)\|_{L^2}<\delta$ and let $u_{g,\delta}$ be the element given by Theorem \ref{appr}. Let $u(t):=U(t,0)u_{g,\delta}$. Then, if $\delta$ is so small that $(M^{1/2}/\kappa)(1+E)\delta^{1/2}<\rho$, the estimate
\begin{equation}
\sup_{t\in[t_*,T]}\|u(t)-\bar u(t)\|_{L^2}\leq K'\exp\left(-N'|\log \delta|^\mu\right)
\end{equation}
holds true (where $K'$ and $N'$ are constants which can be explicitly expressed in terms of $K$, $N$, $\mu$, $\kappa$, $M$ and $E$).\qed
\end{thm}
\begin{rem} The estimate given by Theorem \ref{rec_th_bis} is intermediate between H\umlaut older an logarithmic. Unlike in the Lipschitz case, in the Log-Lipschitz case it is not possible to obtain an approximation rate of H\umlaut older type for $t\geq t_*>0$ (see the remark following Theorem \ref{res_jde_bis}).\end{rem}

In the next section we give two interpretations of this reconstruction method.

\section{Interpretations of the result and relation with similar results} \label{sec:interpr}

\subsection{\it Interpretation as a nonlocal quasi-boundary value problem}

Keeping notations from the previous sections, let $u(t):=U(t,0)u_{g,\delta}$. For $t\in[0,2T]$ define
\begin{equation*}
b_{ij}(t,x):=\begin{cases}a_{ij}(t,x)&\text{for $t\in[0,T]$}\\
a_{ij}(2T-t,x)&\text{for $t\in[T,2T]$}
\end{cases}
\end{equation*}

\begin{equation*}
B(t) u:=-\sum\limits_{j,k=1}^n \partial_{x_j}(b_{jk}(t,x)\partial_{x_k}u)\quad\text{for $t\in[0,2T]$}
\end{equation*}
and
\begin{equation*}
w(t):=\begin{cases}u(t)&\text{for $t\in[0,T]$}\\
V(2T-t,T)u(T)&\text{for $t\in[T,2T]$}
\end{cases}
\end{equation*}
By elliptic regularity and reminding that $V(2T-t,T)=U(T,2T-t)^*$, it is easy to check that $w(t,x)$ solves the (well posed) {\em nonlocal quasi-boundary value problem}
\begin{equation*}
\begin{cases} \partial_t w -\sum\limits_{j,k=1}^n \partial_{x_j}(b_{jk}(t,x)\partial_{x_k}w) = 0&\text{in $[0,2T]$}\\
\alpha(-B(0)w(0)+w(0))+w(2T)=g_2\end{cases}
\end{equation*}
or equivalently
\begin{equation}\label{aux1}
\begin{cases} \partial_t w -\sum\limits_{j,k=1}^n \partial_{x_j}(b_{jk}(t,x)\partial_{x_k}w) = 0&\text{in $[0,2T]$}\\
\alpha(\partial_t w(0)+w(0))+w(2T)=g_2&\end{cases}
\end{equation}
with $\alpha=M\kappa\delta$, where $g_2=U(T,0)^*g=V(0,T)g$. Notice that $g_2=\eta(2T)$, where $\eta(t)$ is the solution of the (well posed) {\em forward} Cauchy problem
\begin{equation}\label{aux2}
\begin{cases} \partial_t\eta -\sum\limits_{j,k=1}^n \partial_{x_j}(b_{jk}(t,x)\partial_{x_k}\eta) = 0&\text{in $[T,2T]$}\\
\eta(T)=g&\end{cases}
\end{equation}
Therefore finding $u_{g,\delta}$ amounts to:
\begin{itemize} 
\item solve the forward Cauchy problem (\ref{aux2});
\item set $g_2:=\eta(2T)$;
\item solve the nonlocal quasi-boundary value problem (\ref{aux1}) with $\alpha=M\kappa\delta$;
\item set $u_{g,\delta}:= w(0)$.
\end{itemize}
Interpreted in this way our approach regains its similarities with the one of \cite{haovanduc}.

\subsection{\it Variational interpretation}
We have seen that $u_{g,\delta}$ is obtained via Lax-Milgram Theorem, and it is the unique minimizer of the functional
\begin{equation}
{\mathcal L}_\alpha (u)=\frac12b_\alpha(u,u)-\langle U(T,0)^*g,u\rangle_{L^2},\quad u\in H^1(\Bbb R^n),
\end{equation}
with $\alpha=M\kappa\delta$. Spelled out in full, the latter reads
\begin{multline*}
{\mathcal L}_\alpha (u)=\frac\alpha2
\left(\int_{\Bbb R^n}\sum\limits_{i,j=1}^n a_{ij}(0,x)\partial_{x_i}u(x)\,\partial_{x_j}u(x)\,dx +\int_{\Bbb R^n}u(x)^2\,dx\right)\\
+\frac12\int_{\Bbb R^n}\left|(U(T,0)u)(x)\right|^2\,dx-\int_{\Bbb R^n}(U(T,0)^*g)(x)u(x)\,dx.
\end{multline*}
Adding the constant $\frac12\|g\|_{L^2}^2$, we obtain that $u_{g,\delta}$ is the unique minimizer of 
\begin{multline*}
\tilde{\mathcal L}_\alpha (u)=\frac\alpha2
\left(\int_{\Bbb R^n}\sum\limits_{i,j=1}^n a_{ij}(0,x)\partial_{x_i}u(x)\,\partial_{x_j}u(x)\,dx +\int_{\Bbb R^n}u(x)^2\,dx\right)\\
+\frac12\int_{\Bbb R^n}\left|(U(T,0)u)(x)-g(x)\right|^2\,dx
\end{multline*}
with $\alpha=M\kappa\delta$. This means, roughly, that $u_{g,\delta}$ is chosen in order to make $\|U(T,0)u-g\|^2+\delta\|u\|_{H^1}^2$ as small as possible.

Interpreted in this way our approach regains its similarities with the one of \cite{liu}.

\subsection{\it Comparison with similar results}
In the Log-Lipschitz case considered in this paper the stability estimates of theorems \ref{res_jde} and \ref{res_jde_bis} are obtained through a rather sophisticated use of the Fourier transform and by mean of paradifferential calculus. This makes the results  not directly exportable to the case of equations on bounded domains, where such tools are not available. That being said, we also point out that though the exponent $\theta$ in (\ref{intro_dis}) and (\ref{intro_rate}) could in principle be tracked through the computations in \cite{CDSP2022}, it is in general small and far from optimal. We discuss briefly how this convergence rate relates to similar results concerning autonomous equations and non-autonomous equations whose coefficients are Lipschitz continuous in $t$. For autonomous parabolic equations on bounded domains, in \cite{denche} Denche and Bessila, with a $H^2$ a-priori bound, obtained approximate solutions with a convergence rate of type (\ref{intro_rate}) with $\theta=1$, while in \cite{HDS} H\`ao,  Duc, and Sahli, with a $H^1$ a-priori bound, obtained a convergence rate of type (\ref{intro_rate}) with
$\theta=1/2$;  along the same line, in \cite{HDL} H\`ao,  Duc, and Lesnic proved several stability estimates up
to 0 of optimal order. The same type of optimal convergence rate was obtained by Liu and Xiao in the paper \cite{liu}, where the authors presented also a very accurate numerical treatment of their approximation result. For the non-autonomous case only few results are available (provided the coefficients are Lipschitz continuous in $t$). As mentioned in the Introduction, in the paper \cite{haovanduc} H\`ao and Duc obtained approximate solutions under a $L^2$ a-priori bound, with a convergence rate of H\umlaut older type for $t\geq t_*>0$, but with no information at $t=0$. In the paper \cite{ImaYama2014} Imanuvilov and Yamamoto proved a conditional stability result of logarithmic type up to $t=0$, which suggests that an approximation scheme with $H^1$ a-priori bound and a convergence rate of type (\ref{intro_rate}) with $\theta=1/4$ can be developed. Finally, in the recent paper \cite{DMA} Duc, Muoi and Anh obtained a convergence rate of type (\ref{intro_rate}) with $\theta=1/4$ under $H^1$ a-priori bound, for a particular non-autonomous parabolic equation. As we stressed above, in our result the exponent $\theta$ is in general small and far from optimal. This is ultimately due to the fact that in the Log-Lipschitz case it is not possible to obtain a conditional stability estimate of H\umlaut older type for $t\geq t_*>0$ (see the remark following Theorem \ref{res_jde_bis}).

\end{document}